# An Accelerated-Decomposition Approach for Security-Constrained Unit Commitment with Corrective Network Reconfiguration

Arun Venkatesh Ramesh, *Student Member, IEEE*, Xingpeng Li, *Member, IEEE* and Kory W. Hedman, *Senior Member, IEEE*

***Abstract*— Security-constrained unit commitment (SCUC) model is used for power system day-ahead scheduling. However, current SCUC model uses a static network to deliver power and meet demand optimally. A dynamic network can provide a lower optimal cost and alleviate network congestion. However, due to the computational complexity and the lack of effective algorithms, network reconfiguration has not been included in the SCUC model yet. This paper presents a novel approach to handle the computational complexity in security-constrained unit commitment (SCUC) with corrective network reconfiguration (CNR) while considering the scalability through accelerated-decomposition approach with fast screening non-critical sub-problems of SCUC-CNR. The proposed approach provides substantial computational benefits and is also applicable to SCUC. Simulation results on the IEEE 24-bus system show that the proposed methods are substantially faster without the loss in solution quality while the scalability benefits are demonstrated using larger cases: the IEEE 73-bus system, IEEE 118-bus system and Polish system.**

***Index Terms*— Accelerated-decomposition approach, Benders decomposition, Corrective transmission switching, Flexible transmission, Mixed-integer linear programming, Network reconfiguration, Post-contingency congestion relief, Security-constrained unit commitment, Topology control.**

## NOMENCLATURE

Indices:

| | |
|---|---|
| $c$ | Line contingency index. |
| $g$ | Generator index. |
| $k$ | Transmission element (line or transformer) index. |
| $n$ | Bus index. |
| $N(g)$ | Bus location of generator g. |
| $t$ | Time period index. |

Sets:

| | |
|---|---|
| $C$ | Set of non-radial transmission contingencies. |
| $G$ | Set of generators. |
| $g(n)$ | Set of generators connecting bus $n$. |
| $K$ | Set of all transmission element. |
| $K_r$ | Set of reconfigurable non-radial lines. |
| $N$ | Set of all buses. |
| $T$ | Set of Time intervals. |
| $\delta^+(n)$ | Set of lines with bus $n$ as receiving bus. |
| $\delta^-(n)$ | Set of lines with bus $n$ as sending bus. |
| $\Omega^{cri}$ | Set of all critical sub-problems. |
| $\Omega_1^{inf}$ | Set of infeasible PCFC sub-problems. |
| $\Omega_2^{inf}$ | Set of infeasible NR-PCFC sub-problems. |
| $\psi$ | Cut-set determined for all sub-problems. |

Parameters:

| | |
|---|---|
| $b_k$ | Susceptance of line $k$. |
| $c_g$ | Linear cost for generator $g$. |
| $c_g^{NL}$ | No-load cost for generator $g$. |
| $c_g^{SU}$ | Start-up cost for generator $g$. |
| $d_{n,t}$ | Predicted demand of bus $n$ in time period $t$. |
| $DT_g$ | Minimum down time for generator $g$. |
| $M$ | A big real number. |
| $P_g^{max}$ | Maximum capacity of generator $g$. |
| $P_g^{min}$ | Minimum capacity of generator $g$. |
| $P_k^{emax}$ | Emergency thermal line limit for line $k$. |
| $P_k^{max}$ | Long-term thermal line limit for line $k$. |
| $P_{g,t}^{MUC}$ | Output of generator $g$ in period t obtained from MUC. |
| $P_{k,t}^{MUC}$ | Flow in line $k$ in time period $t$ obtained from MUC. |
| $R_g^{10}$ | 10-minute outage ramping limit of generator $g$. |
| $R_g^{hr}$ | Regular hourly ramping limit of generator $g$. |
| $R_g^{SD}$ | Shut-down ramping limit of generator $g$. |
| $R_g^{SU}$ | Start-up ramping limit of generator $g$. |
| $UT_g$ | Minimum up time for generator $g$. |
| $u_{g,t}^{MUC}$ | Generator $g$ status in period $t$ obtained from MUC. |
| $v_{g,t}^{MUC}$ | Generator $g$ start-up in period $t$ obtained from MUC. |
| $Z_{max}$ | Maximum number of CNR actions per sub-problem. |
| $\theta_{ref,t}$ | Phase angle of reference bus in time period $t$. |

Variables:

| | |
|---|---|
| $F_{k,c,t}^+, F_{k,c,t}^-$ | Dual variables of line $k$'s contingent max and min limit constraints for contingency $c$ and period $t$. |
| $P_{g,c,t}$ | Output of generator $g$ in period $t$ after outage of line $c$ |
| $P_{k,c,t}$ | Flow in line $k$ in period $t$ after outage of line $c$. |
| $P_{g,t}$ | Output of generator $g$ in time period $t$. |
| $P_{k,t}$ | Lineflow of line $k$ in time period $t$. |
| $r_{g,t}$ | Reserve from generator $g$ in time period $t$. |
| $S_{k,c,t}$ | Dual variable of line $k$'s contingent power flow constraint for contingency $c$ and period $t$. |
| $u_{g,t}$ | Commitment status of generator $g$ in time period $t$. |
| $v_{g,t}$ | Start-up variable of generator $g$ in time period $t$. |
| $z_{c,t}^k$ | Line status for line $k$ in period $t$ after outage of line c. |
| $\alpha_{g,c,t}^+, \alpha_{g,c,t}^-$ | Dual variables of generator $g$ contingent max and min capacity constraint, for contingency $c$ and period $t$. |
| $\beta_{g,c,t}^+, \beta_{g,c,t}^-$ | Dual variables of generator $g$ contingent reserve max and min constraint, for contingency $c$ and period $t$. |
| $\lambda_{n,c,t}$ | Dual variable of bus $n$'s power balance constraint for contingency $c$ and period $t$. |
| $\theta_{m,c,t}$ | Phase angle of bus $m$ in period $t$ after outage of line $c$. |
| $\theta_{n,c,t}$ | Phase angle of bus $n$ in period $t$ after outage of line $c$. |
| $\theta_{m,t}$ | Phase angle of bus $m$ in time period $t$. |
| $\theta_{n,t}$ | Phase angle of bus $n$ in time period $t$. |

Arun Venkatesh Ramesh and Xingpeng Li are with the Department of Electrical and Computer Engineering, University of Houston, Houston, TX, 77204, USA. Kory W. Hedman is with the School of Electrical, Computer and Energy Engineering, Arizona State University, Tempe, AZ, 85287, USA (e-mail: aramesh4@uh.edu; xingpeng.li@asu.edu; kwh@myuw.net).

## I. Introduction

The electric power needs to be generated, transferred and utilized concurrently. This requires state of the art approaches that optimize the scheduling before-hand to ensure reliable power supply, save cost and avoid resource wastage. This stresses on smarter algorithms to effectively utilize the flexibility in the power system that includes the network.

The transmission network is built with a lot of redundancy since it generally considers future demand growth and meets high reliability standards. Though the transmission flexibility through network reconfiguration (NR) is a cheap and quick action, the transmission elements are mostly considered as a static element in day-ahead scheduling especially for congestion management [1]. Apart from system reliability, NR provides significant cost-saving benefits [1]-[2] and network congestion alleviation by rerouting the network flows [3].

The North American Electric Reliability Corporations (NERC) sets *N*-1 standards to ensure stability of the system which are capable of handling typical system disturbances [4]. Here, system operators utilize both preventive and corrective actions to handle the above uncertainties. A preventive NR action is implemented prior to the contingency to avoid line flows exceeding emergency ratings following a contingency. This means the system is set to handle the contingency without any control actions such as re-dispatch following a contingency.

A corrective NR action is implemented after the disturbance has occurred to move the system from emergency state to normal-secure or normal-insecure state. NR when used as a corrective action can re-route the line flows and relieve post-contingency network congestion, which may allow cheaper generators to produce more power. It can be noted that corrective network reconfiguration (CNR) is only utilized when a contingency occurs, and the associated post-contingency network is overloaded. Also, CNR avoids large system disturbances that NR leads to and reduces circuit breaker degradation [5].

The importance of CNR is seen through several industrial examples based on historical or simulated control schemes. Presently, transmission operators follow the procedure for relieving network congestion based on experience such as [6]-[7] for PJM and [8] for ISO New England rather than sound systematic methods especially during contingencies. Such actions were used during disasters like Superstorm Sandy [9].

The consideration of CNR in SCUC was first seen in [10], however, the approach considers only a small subset of reconfigurable candidate lines and contingencies. This was bettered in [10] which considers all transmission contingencies all non-radial lines as reconfigurations candidate lines. However, this extensive approach is not scalable to large systems. Prior research [11] demonstrates that a co-optimized CNR method leads to significant cost saving and network congestion alleviation and [12] shows that CNR can benefit integration of renewable energy. Moreover, the transmission line overload reduction and market surplus benefits were realized effectively through CNR in [13]-[15].

Network flexibility can be introduced in both real-time and day-ahead operations in the bulk power system. Due to the complexity, it can be noted from [16]-[18], NR is implemented by various heuristic methods to obtain quick results. [19] utilizes three concurrent NR actions to improve performance. In real-time scenarios, [5] presents a framework for integrating CNR with real-time contingency analysis and [20]-[21] proposed an enhanced energy management system with inclusion of a CNR module that can seamlessly and practically connect with real-time contingency analysis and security-constrained economic dispatch.

In day-ahead scenario, the security-constrained unit commitment (SCUC) is run to obtain an economical viable solution along with the day-ahead generator commitment and dispatch schedule. Since SCUC is used in both regulated and deregulated environments, the algorithm developed in the paper can be implemented in either business environment.

One main reason for not including NR/CNR is the increase in complexity of the *N*-1 SCUC model as it introduces additional binary variables to the mixed integer linear programming (MILP) problem. Here, decomposing the SCUC by iterative multi-stage approaches or using heuristic techniques is beneficial for algorithm performance. [22] recognize NR benefits with a small subset of reconfigurable assets. [23] proposes a co-optimized method which enhances *N*-1 security by considering both a preventive optimal NR scheduling and a CNR rescheduling that tolerate short-term overloads in post-contingency scenarios. [24]-[25] detail a two-stage SCUC with NR that can be solved iteratively for large-scale power systems. [26] proposes an iterative fast SCUC method to compute for each hour and provide the resulting solution as a starting point for the original SCUC. Typically decomposing the SCUC is implemented using Benders decomposition algorithm (BDA) or column and constraint generation algorithms (CCGA).

BDA can effectively reduce the complexity of SCUC by decomposing it as a master problem and associated sub-problems. [27] solves a stochastic-SCUC problem which implements NR to mitigate wind uncertainty and considers an AC optimal power flow through linearized network losses by utilizing BDA to reduce the problem complexity. [28] implements a multi-stage discrete approach through BDA acceleration techniques to include emerging technologies in SCUC. However, [26]-[28] does not consider NR/CNR. In [29], a sequential extensive approach to implement CNR in *N*-1-1 SCUC is considered, which is not scalable.

Several authors have indicated slow convergence and sought for heuristics to speed up BDA. In [30], a mixed-integer non-linear problem of long-term planning of distributed generation and reconfiguration of distribution system was solved by using BDA and accelerated by considering both feasibility and optimality cuts. In [31], ordered sets, pre-solve and warm-start techniques were used for short-term hydropower maintenance scheduling. However [30]-[31] do not perform any day-ahead market operations. In [32], a security-constrained optimal power flow was solved using heuristics and parallel solving of sub-problems using BDA. In [33], the concept of strong cuts using sensitivity factors were introduced to reduce iterations for solving SCUC. However, NR/CNR was never considered in [30]-[33].

CCGA, like BDA, also helps in reducing the problem complexity in SCUC. It is seen from [34], that energy and reserves are co-optimized while considering CNR and pre-contingency NR (PNR). A multi-level nested CCGA is used in [34]. Though this method addresses several extensions and considers umbrella contingencies, only one worst-case critical contingency per iteration is addressed to provide convergence for an exact solution and not compromise solve time. Therefore, the solution obtained from this method does not provide reconfiguration or re-dispatch solution for all contingencies. If more contingencies are required to be addressed then the method proposed in [34] requires further iterations which may lead to substantial increase in solve time.

The absence of a reliable and scalable algorithm which implements CNR for large-scale practical power systems with a capability of considering all system contingencies is the research gap addressed in this paper. This is implemented without compromising solve time while also providing the solution for re-dispatch and CNR (if present) for all contingencies. This paper proposes decomposition approaches for SCUC and SCUC-CNR; and the proposed accelerators and screener can enhance the proposed methods. The contributions of this work are presented as follows:

- The proposed accelerated-decomposition approaches to SCUC and SCUC-CNR are scalable to large-scale power systems such as Polish network. In particular, the non-convexity with the addition of CNR variables through Benders decomposition is innovatively addressed while maintaining solution quality.
- The developed contingency screener is fast and considers the entire list of non-radial lines as the contingency list, and it efficiently identifies critical contingencies that may lead to overloads.
- The ranked closest branches to contingency element priority list can obtain optimal CNR actions quickly while ensuring feasibility of post-contingent scenarios.
- The proposed decomposition approaches provide the benefit of easy integration of additional accelerators, technologies and techniques as pointed in Section VIII.

The rest of this paper is organized as follows. Section II presents the extensive formulation. Section III describes the decomposition of the extensive formulation. Section IV models the proposed approaches. Section V details the testing data and Section VI analyzes the results. Finally, Section VII concludes the paper and the potential extensions to this work are presented in Section VIII.

## II. Extensive Formulation

The objective of SCUC is to minimize operational cost of generators (1). This is accomplished subject to base-case and post-contingency constraints which are co-optimized in an extensive formulation. The base-case generation constraints are modelled in (2)-(11). Constraints (2) and (3) represent the minimum and maximum generation limits; (4) and (5) enforce the reserve requirements; (6) and (7) are the hourly ramping limits; (8) and (9) are the generator min-up and min-down time constraints. Generator start-up indication variable is defined in (10). The generator commitment and start-up variables are bound by binary integrality constraints as shown in (11). The base-case physical power flow constraint is represented through (12)-(14). (12) depicts the power flow calculation; (13) represents the long-term thermal limits of transmission elements; and (14) enforces nodal power balance. Slack equation, (15), is added to define the reference phase angle in the base-case solution.

*Objective:*

$$Min\ \sum_g \sum_t \left( c_g P_{g,t} + c_g^{NL} u_{g,t} + c_g^{SU} v_{g,t} \right) \tag{1}$$

*s.t.:*

*Base case modeling of generation:*

$$P_g^{min} u_{g,t} \leq P_{g,t}, \forall g, t \tag{2}$$

$$P_{g,t} + r_{g,t} \leq P_g^{max} u_{g,t}, \forall g, t \tag{3}$$

$$0 \leq r_{g,t} \leq R_g^{10} u_{g,t}, \forall g, t \tag{4}$$

$$\sum_{q \in G} r_{q,t} \geq P_{g,t} + r_{g,t}, \forall g, t \tag{5}$$

$$P_{g,t} - P_{g,t-1} \leq R_g^{hr} u_{g,t-1} + R_g^{SU} v_{g,t}, \forall g, t \tag{6}$$

$$P_{g,t-1} - P_{g,t} \leq R_g^{hr} u_{g,t} + R_g^{SD} (v_{g,t} - u_{g,t} + u_{g,t-1}), \forall g, t \tag{7}$$

$$\sum_{q=t-UT_g+1}^{t} v_{g,q} \leq u_{g,t}, \forall g, t \geq UT_g \tag{8}$$

$$\sum_{q=t+1}^{t+DT_g} v_{g,q} \leq 1 - u_{g,t}, \forall g, t \leq T - DT_g \tag{9}$$

$$v_{g,t} \geq u_{g,t} - u_{g,t-1}, \forall g, t \tag{10}$$

$$v_{g,t}, u_{g,t} \in \{0,1\}, \forall g, t \tag{11}$$

*Base case modeling of power flow:*

$$P_{k,t} - b_k \left( \theta_{n,t} - \theta_{m,t} \right) = 0, \forall k, t \tag{12}$$

$$-P_k^{max} \leq P_{k,t} \leq P_k^{max}, \forall k, t \tag{13}$$

$$\sum_{g \in g(n)} P_{g,t} + \sum_{k \in \delta^+(n)} P_{k,t} - \sum_{k \in \delta^-(n)} P_{k,t} = d_{n,t}, \forall n, t \tag{14}$$

$$\theta_{ref,t} = 0\ \forall t \tag{15}$$

The post-contingency constraints are modelled through the post-contingency generation constraints, (16)-(19), and post-contingency power flow constraints for non-radial lines, (20)-(22). The post-contingency constraints ensure that the disturbance in the system is handled in 10 minutes. Here, (16)-(17) represent the 10-minute ramp up/down limit; (18)-(19) model the minimum and maximum limits of the generator. (20) calculates the post-contingency line flows; (21) enforces the emergency rating of the transmission element. Finally, (22) represents the nodal power balance in the post-contingency case.

*Post-contingency 10-minute ramping restriction on generation and modeling of contingencies:*

$$P_{g,t} - P_{g,c,t} \leq R_g^{10} u_{g,t}, \forall g, c \in C, t \tag{16}$$

$$P_{g,c,t} - P_{g,t} \leq R_g^{10} u_{g,t}, \forall g, c \in C, t \tag{17}$$

$$P_g^{min} u_{g,t} \leq P_{g,c,t}, \forall g, c \in C, t \tag{18}$$

$$P_{g,c,t} \leq P_g^{max} u_{g,t}, \forall g, c \in C, t \tag{19}$$

*Post-contingency modeling of power flow:*

$$P_{k,c,t} - b_k \left( \theta_{n,c,t} - \theta_{m,c,t} \right) = 0, \forall k, c \in C, t \tag{20}$$

$$-P_k^{emax} \leq P_{k,c,t} \leq P_k^{emax}, \forall k, c \in C, t \tag{21}$$

$$\sum_{g \in g(n)} P_{g,c,t} + \sum_{k \in \delta^+(n)} P_{k,c,t} - \sum_{k \in \delta^-(n)} P_{k,c,t} = d_{n,t}, \forall n, c \in C, t \tag{22}$$

If CNR is modelled, then the post-contingency transmission constraints, (20)-(21), are replaced with the post-contingency line flow equations and limits with CNR, (23)-(26). The linearity of post-contingency power flow equations, (23)-(24), are maintained with the *big-M* method. The binary decision variable, $z_{c,t}^{k}$, represents the CNR action where a value of 0 represents line is disconnected from the system and a value of 1 indicates line is available. These contingencies are modelled for all non-radial lines. A restriction on the number of CNR actions in each post-contingency case is introduced through (26) to reduce system disturbance.

*Post-contingency modeling of power flow when CNR is incorporated:*

$$P_{k,c,t} - b_k\left(\theta_{n,c,t} - \theta_{m,c,t}\right) + \left(1 - z_{c,t}^{k}\right)M \geq 0, \forall k, c \in C, t \tag{23}$$

$$P_{k,c,t} - b_k\left(\theta_{n,c,t} - \theta_{m,c,t}\right) - \left(1 - z_{c,t}^{k}\right)M \leq 0, \forall k, c \in C, t \tag{24}$$

$$-P_k^{emax} z_{c,t}^{k} \leq P_{k,c,t} \leq z_{c,t}^{k} P_k^{emax}, \forall k, c \in C, t \tag{25}$$

$$\textstyle\sum_k \left(1 - z_{c,t}^{k}\right) \leq Z_{max}, \forall k, c \in C, t \tag{26}$$

Based on the above constraints, the SCUC formulation is represented by (1)-(22) and SCUC-CNR is represented by (1)-(19) and (22)-(26).

## III. Decomposition of Extensive Formulation

The extensive formulations for SCUC and SCUC-CNR are decomposed as master problem and several sub-problems to reduce the computational burden using BDA approach. The resulting master problem is a relaxed MILP obtained by the base-case constraints which provides the dispatch and commitment schedules and is referred to as master unit commitment (MUC) problem in this paper. The sub-problems are LP problems obtained from the post-contingency model of SCUC which checks the feasibility of base-case solution are termed as post-contingency feasibility check (PCFC). The BDA approach is a row-generalized approach, meaning extra cut constraints are added to the master problem to link the sub-problems using duality theory. In SCUC, the cut constraint which is added to the relaxed master problem is aggregated from the dual values of the PCFC sub-problems that violates the physical constraints at the end of each iteration. Subsection III.A details the Master problem and subsection III.C details the PCFC sub-problem.

In comparison with SCUC, SCUC-CNR varies in the post-contingency power flow constraints, and the introduction of reconfiguring variable, $z_{c,t}^{k}$, makes the prior derivation impossible as it leads to an MILP-based sub-problem. This is overcome by adding an extra sub-problem, network-reconfigured PCFC (NR-PCFC), which iterates through reconfiguration action, one at a time, to make the problem an LP as presented in sub-section III.D. Although this is a heuristic method rather than an exact algorithm that guarantees optimality, it is a very effective method as demonstrated in section VI.

Two accelerators were identified to make the typical-decomposition approach solve faster: ranked closest branches to contingency element (CBCE) list [5], and critical subproblem identification. Firstly, it was observed that only a subset of the contingency sub-problems are critical, and an accelerator was developed to identify critical sub-problems to reduce computational burden. This accelerator can be implemented for both SCUC and SCUC-CNR and is represented in sub-section III.B as critical sub-problem screener (CSPS). Secondly, the CNR actions can be implemented through the CBCE list, a ranked priority list of 20 closest branches to each contingent element in the network to obtain quick feasible results for CNR. The CBCE list is only used in NR-PCFC sub-problems and therefore, it is only used in the proposed methods implementing SCUC-CNR for large networks.

### A. Master Unit Commitment

The MUC problem is represented through (1)-(15) and (27). It obtains the base-case solution which provides the generator commitment and dispatch for all periods that are then used in sub-problems. The MUC is an MILP problem.

$$\begin{aligned} &\textstyle\sum_{g\in G}\big(P_g^{max} u_{g,t}^{fix} \alpha_{g,c,t}^{+} - P_g^{min} u_{g,t}^{fix} \alpha_{g,c,t}^{-} + \\ &(R_g^{10} u_{g,t}^{fix} - P_{g,t}^{fix})\beta_{g,c,t}^{+} + (R_g^{10} u_{g,t}^{fix} + \\ &P_{g,t}^{fix})\beta_{g,c,t}^{-}\big) + \textstyle\sum_{k\in K/\{c\}}\big(P_k^{emax}(F_{k,c,t}^{+} + F_{k,c,t}^{-}) + \\ &0(S_{k,c,t})\big) + \textstyle\sum_{n\in N} d_{n,t}\lambda_{n,c,t} \leq 0, \forall \psi \end{aligned} \tag{27}$$

Based on the BDA approach, (27) represents the cuts associated with infeasible sub-problems from post-contingency feasibility check (PCFC) using duality theory. In the proposed methodologies, the cuts are purely created with the dual-variables of the sub-problem PCFC (29)-(38). The NR-PCFC sub-problem only verifies feasibility of the sub-problems that are infeasible by PCFC through various re-configurations; it was observed that this implementation has little impact on the solution optimality.

The cut-set, $\psi$, is obtained after PCFC in the case of SCUC or after NR-PCFC in the case of SCUC-CNR for each iteration. Once MUC is solved, the set $\Omega^{cri}$ is initialized with a complete list of sub-problem $c \in C$, $t \in T$. The MUC commitment and dispatch are passed on to each sub-problem in CSPS, PCFC and NR-PCFC.

### B. Critical Sub-Problem Screener

The purpose of the critical sub-problem screener is to quickly screen out non-critical sub-problems before PCFC and NR-PCFC. Fig. 1 (a) depicts the flow of CSPS.

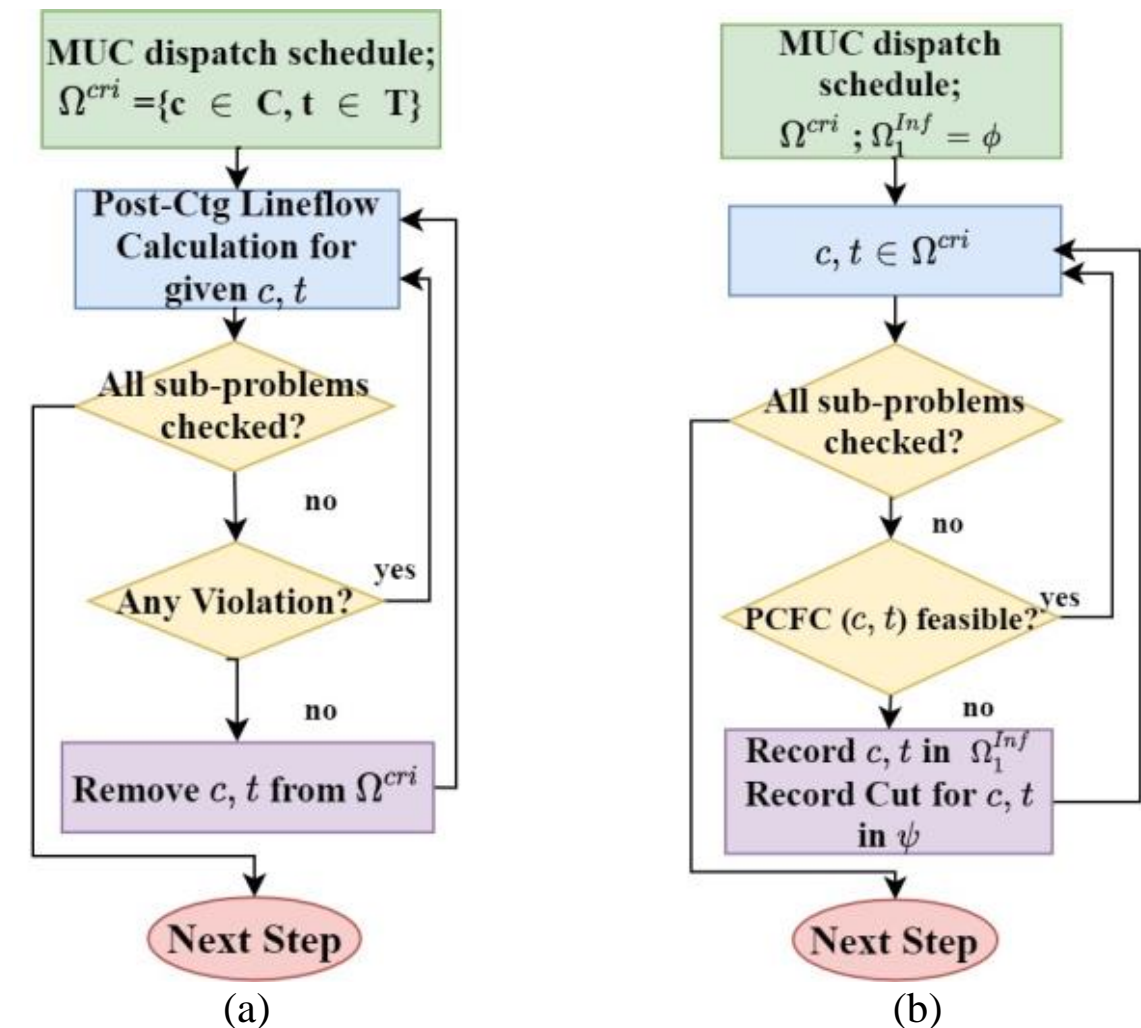

Fig. 1. Flowchart: (a) Critical sub-problem screener (b) Post-contingency feasibility check.

Post-contingent line flows for each sub-problem in critical set, $\Omega^{cri}$, are obtained through the predetermined line outage distribution factor (LODF), (28). The contingent line flows are then compared against the emergency line limit for violations.

The non-critical sub-problems determined by CSPS are removed from the set $\Omega^{cri}$ leaving only critical sub-problems. CSPS is very fast since it only involves a limited number of simple algebraic calculations rather than complex optimization.

$$P_{k,c,t} = P_{k,t}^{MUC} + LODF_{k,c}\left(P_{c,t}^{MUC}\right), \forall k \tag{28}$$

### C. Post-Contingency Feasibility Check

The sub-problem PCFC is represented by (29)-(38), which is derived from (16)-(22). The goal of PCFC is to check system feasibility for sub-problems in set $\Omega^{cri}$ by conducting emergency generation re-dispatch without CNR. This is done by minimizing the non-negative slack variable, $s_1$, which indicates the feasibility of the sub-problem. If $s_1$ is exactly zero, then the problem is feasible; otherwise it is infeasible. If PCFC fails feasibility, the respective sub-problem ($c$, $t$) will be recorded to set $\Omega_1^{inf}$ along with respective cut in the cut-set, $\psi$. Fig. 1 (b) depicts the flow of PCFC.

*Objective:*

$$Min\ s_1 \tag{29}$$

*s.t.:*

*Post-contingency generation modeling for a given contingency c in time period t in set $\Omega^{cri}$:*

$$-P_{g,c,t} + s_1\left(R_g^{10}u_{g,t}^{MUC} - P_{g,t}^{MUC}\right) \le R_g^{10}u_{g,t}^{MUC} - P_{g,t}^{MUC}, \forall g \qquad (\beta_{g,c,t}^{-}) \tag{30}$$

$$P_{g,c,t} + s_1\left(R_g^{10}u_{g,t}^{MUC} + P_{g,t}^{MUC}\right) \le R_g^{10}u_{g,t}^{MUC} + P_{g,t}^{MUC}, \forall g \qquad (\beta_{g,c,t}^{+}) \tag{31}$$

$$P_g^{min}u_{g,t}^{MUC} \le P_{g,c,t} + s_1\left(P_g^{min}u_{g,t}^{MUC}\right), \forall g \qquad (\alpha_{g,c,t}^{-}) \tag{32}$$

$$P_{g,c,t} + s_1\left(P_g^{max}u_{g,t}^{MUC}\right) \le P_g^{max}u_{g,t}^{MUC}, \forall g \qquad (\alpha_{g,c,t}^{+}) \tag{33}$$

*Post-contingency modeling of power flow for a given contingency c in time period t in set $\Omega^{cri}$:*

$$P_{k,c,t} - b_k\left(\theta_{n,c,t} - \theta_{m,c,t}\right) = 0, \forall k \in K/\{c\} \qquad (S_{k,c,t}) \tag{34}$$

$$P_{c,c,t} = 0 \tag{35}$$

$$-P_k^{emax} \le P_{k,c,t} - s_1(P_k^{emax}), \forall k \qquad (F_{k,c,t}^{-}) \tag{36}$$

$$P_{k,c,t} + s_1(P_k^{emax}) \le P_k^{emax}, \forall k \qquad (F_{k,c,t}^{+}) \tag{37}$$

$$\sum_{g\in g(n)} P_{g,c,t} + \sum_{k\in\delta^{+}(n)} P_{k,c,t} - \sum_{k\in\delta^{-}(n)} P_{k,c,t} + s_1(d_{n,t}) = d_{n,t}, \forall n \qquad (\lambda_{n,c,t}) \tag{38}$$

### D. Network-Reconfigured Post-Contingency Feasibility Check

The sub-problem NR-PCFC that includes CNR is represented by (29)-(33), (35)-(38) and (39)-(40). The goal of network-reconfigured post-contingency feasibility check is to check system feasibility with CNR for the set $\Omega_1^{inf}$. The feasibility is checked by switching one non-radial transmission element $j$ at a time from the network. As mentioned previously, the non-convexity of reconfiguring decision variable, $z_{c,t}^{k}$, is overcome by iterating through CBCE list or complete enumeration of reconfigurable non-radial lines set $K_r$ (benchmark to CBCE) one at a time to identify the reconfigured network topology.

For each scenario (with line $j$ removed from the network) of NR-PCFC, the non-negative slack variable $s_1$ indicates the sub-problem feasibility similar to PCFC. If the sub-problem is feasible for one such scenario, then the sub-problem ($c$,$t$) is feasible through CNR and is removed from the cut-set, $\psi$, obtained from PCFC. Record the line $j$ selected from the CBCE list that facilitates CNR. If no switching scenario leads to a feasible solution for sub-problem ($c$,$t$), then the infeasible sub-problem will be recorded in set $\Omega_2^{inf}$. Fig. 2 depicts the flow of NR-PCFC.

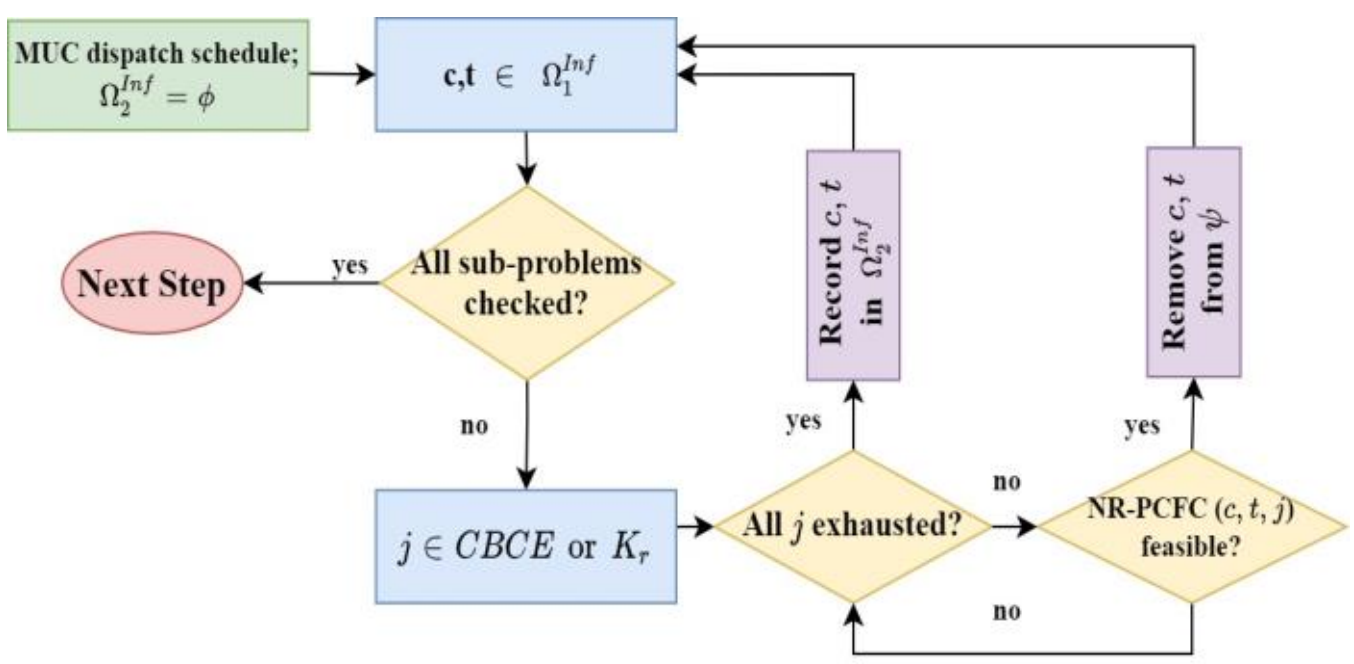

Fig. 2. Network-reconfigured post-contingency feasibility check.

*Post-contingency model of power flow for a given contingency c in time period t in set $\Omega_1^{inf}$ and line j from CBCE:*

$$P_{k,c,t} - b_k\left(\theta_{n,c,t} - \theta_{m,c,t}\right) = 0, \forall k \in K/\{c, j\} \tag{39}$$

$$P_{j,c,t} = 0 \tag{40}$$

## IV. Proposed Methods

This paper compares extensive formulations against decomposition approaches of SCUC and SCUC-CNR and details the benefits of decomposed approaches. The extensive formulations were discussed in Section II. The decomposition approaches to SCUC and SCUC-CNR are explained in the following sub-sections. This paper proposes two decomposed approaches for SCUC namely: typical-decomposition approach to SCUC (T-SCUC) and accelerated-decomposition approach to SCUC (A-SCUC). Along with the above proposed methods, this paper also proposes two decomposed approaches for SCUC-CNR which perform network reconfiguration as a corrective action namely: typical-decomposition approach to SCUC-CNR (T-SCUC-CNR) and accelerated-decomposition approach to SCUC-CNR (A-SCUC-CNR). The proposed approaches are explained through the decomposed-features of master and sub-problems explained in sub-sections III.A-III.D.

### A. Typical-Decomposition Approach

The proposed typical-decomposition approach by using MUC and PCFC only for SCUC whereas SCUC-CNR also utilizes NR-PCFC. The MUC problem is initially solved to obtain the generator commitment and base-case output. The feasibility of each sub-problem in set $\Omega^{cri}$ is checked by post-contingency generation redispatch implemented by PCFC. For the typical-decomposition approach, the set $\Omega^{cri}$ holds the complete list of all sub-problems and the set $\Omega_1^{inf}$ is an empty set at the beginning of each iteration. When the feasibility of a sub-problem is not achieved, it is recorded in the set $\Omega_1^{inf}$.

For SCUC, once all sub-problems are examined, an iteration is completed. The infeasible sub-problem are recorded in set $\Omega_1^{inf}$ at the end of each iteration. The problem is converged when set $\Omega_1^{inf}$ is empty at the end of an iteration.

For SCUC-CNR, the set $\Omega_1^{inf}$ is passed on to NR-PCFC and feasibility of each sub-problem is examined with CNR. If the sub-problem is infeasible then it is recorded in set $\Omega_2^{inf}$. Once all sub-problems in set $\Omega_1^{inf}$ are checked, an iteration is completed.

The respective cuts for infeasible sub-problems in set $\Omega_1^{inf}$ for SCUC and infeasible sub-problems in set $\Omega_2^{inf}$ are formed using the dual value of (19)-(28) and added as (27), respectively, after each iteration.

### *B. Accelerated-Decomposition Approach*

Accelerated-decomposition approach uses MUC, CSPS and PCFC only for SCUC whereas SCUC-CNR also utilizes NR-PCFC. The flow of this approach is similar to typical-decomposition approach, but it is substantially sped through the CSPS, an accelerator to reduce the computational burden by identifying critical sub-problems. The MUC problem is initially solved to obtain the generator commitment and base-case output. With the MUC schedule, the critical sub-problems are identified and recorded in set $\Omega^{cri}$ by using CSPS. Only the critical sub-problems, rather than all sub-problems, are then checked by post-contingency generation redispatch through PCFC.

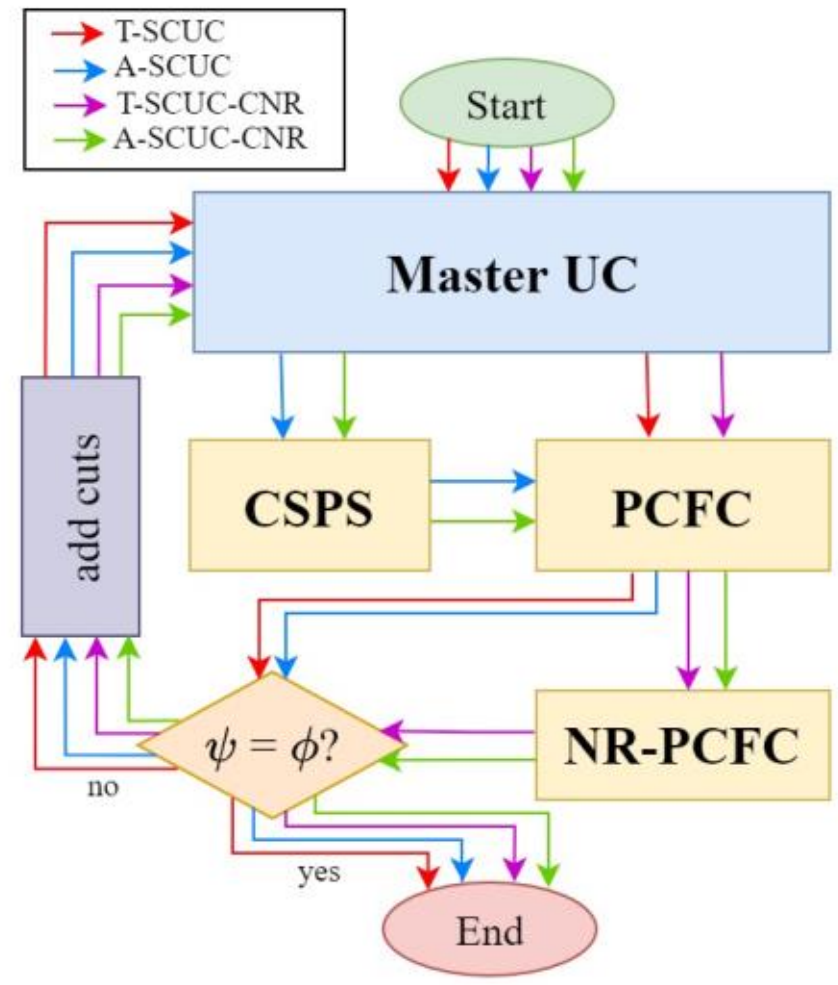

Fig. 3. Flowchart of typical-decomposition approach to SCUC.

**Algorithm 1** Accelerated-decomposition approach to SCUC-CNR

1: Solve MUC and obtain the commitment and dispatch
2: **repeat**
3: cut = ∅;
4: **for** all $t \in T$ **do**
5: **for** all $c \in C$ **do**
6: solve CSPS($c$,$t$)
7: **for** all $k \in K$ **do**
8: **if** $P_{k,c,t}$ violation **then**
9: record ($c$,$t$) in set $\Omega^{cri}$
10: **end if**
11: **end for**
12: **for** all ($c$,$t$) $\in \Omega^{cri}$ **do**
13: solve PCFC($c$,$t$)
14: **if** PCFC($c$,$t$) is infeasible **then**
15: Flag = false
16: **for** line $j$ ∈ CBCE **do**
17: Remove line $j$ from topology
18: solve NR-PCFC($j$,$c$,$t$)
19: **if** NR-PCFC($j$,$c$,$t$) feasible **then**
20: Flag = true; **break**
21: **end if**
22: **end for**
23: **if** Flag = false **then**
24: {cut} = {cut} + {cut of PCFC($c$,$t$)}
25: **end if**
26: **end if**
27: **end for**
28: **end for**
29: **end for**
30: **if** {cut} != ∅ **then**
31: add cut to MUC; solve updated MUC
32: **else**
33: problem converged; report results; **break**
34: **end if**
35: **until** converged

Fig. 3 represents the flow of the proposed typical decomposition approach and proposed accelerated-decomposition approach to SCUC/SCUC-CNR and the pseudo-code is represented in Algorithm 1. Here (i) A-SCUC-CNR is implemented using lines 1-35; (ii) T-SCUC-CNR is implemented through lines 1-5 and 12-35; (iii) T-SCUC is implemented by lines 1-5, 13-14, 24, and 26-35; and (iv) A-SCUC is implemented through lines 1-15 and 23-35.

## V. Test Case Description

The proposed methods, T-SCUC/A-SCUC and T-SCUC-CNR/A-SCUC-CNR were validated against the extensive formulation detailed in SCUC and SCUC-CNR, respectively, on the IEEE 24-bus system with 33 generators and 38 branches [35]. The network includes a total generation capacity of 3,393 MW and the system peak load is 2,265 MW. Furthermore, the IEEE 73-bus system and the Polish system were utilized to show the effectiveness and scalability of T/A-SCUC-CNR. Table I summarizes the test systems.

The IEEE 73-bus system consists of 99 generators and 117 branches [35]. The total generation capacity is 10,215 MW and the system peak load is 8,550 MW. The Polish system, modified to include default min-up/min-down times and ramp-up/ramp-down limits, is used for demonstrating the scalability of the algorithm. It is the largest system used for this work and it consists of 2,383 buses, 327 generators and 2,895 branches [36]. The total generation capacity is 30,053 MW serving a system peak load of 21,538 MW. Two cases of the Polish system, covering a single-hour period and a 24-hour period respectively, are considered. The single-hour period case is effective to compare performance against smaller systems whereas the scalability is shown through the 24-hour period case. For the purpose of demonstrating CNR, only non-radial transmission line contingencies are considered in the *N*-1 SCUC formulation since contingency of radial lines will lead to islanding and system separation; this is consistent with industrial practice. Similarly, CNR actions, at most one action per contingency, considers only non-radial lines as possible reconfiguration actions for the same reason.

In addition, the IEEE 118-bus system with 54 generators and 186 branches, [36], was also considered to draw comparison with state-of-the-art methods to show the efficacy of proposed methods. However, the IEEE 118-bus system lacks the generator data and thermal limits and thus it was modified to include such information.

TABLE I. TEST SYSTEM SUMMARY

| System | Pgen (GW) | Pload (GW) | # bus | #gen | # branch | # radial branch |
|---|---|---|---|---|---|---|
| IEEE 24 | ~3.4 | ~2.1 | 24 | 33 | 38 | 1 |
| IEEE 73 | ~10.2 | ~8.6 | 73 | 99 | 117 | 2 |
| IEEE 118 | ~5.8 | ~3.1 | 118 | 54 | 186 | 7 |
| Polish | ~30.1 | ~21.5 | 2,383 | 327 | 2,895 | 644 |

## VI. RESULTS AND ANALYSIS

The mathematical model is implemented using AMPL and solved using Gurobi. The proposed methods were initially validated, following which sensitivity analysis, scalability and market impact are discussed.

The typical-decomposition and accelerated-decomposition approaches are proposed to address the computational complexity of both SCUC and SCUC-CNR. Especially, the proposed accelerated-decomposition approach by utilizing critical sub-problem screener (CSPS) substantially reduces the solve time and outperforms the typical-decomposition approach and extensive formulation. It can effectively be scalable to handle large power systems. The proposed methods can be easily integrated into the existing practices of SCUC without the loss of solution quality.

Here, SCUC serves as a benchmark for proposed T-SCUC and A-SCUC whereas SCUC-CNR serves as a benchmark for proposed T-SCUC-CNR and A-SCUC-CNR.

### *A. Validation for Proposed Methodologies*

Since the proposed methodologies are all iterative in nature, an accuracy validation was performed to test the robustness against non-iterative extensive formulations. A MIPGAP of 0.00 was utilized on the congested network of IEEE 24-bus system for 24-hour period and the SCUC results are tabulated in Table II and SCUC-CNR results are tabulated in Table III. It was observed from Table II that the results for SCUC, T-SCUC and A-SCUC are the same. Similarly, the solutions obtained from SCUC-CNR, T-SCUC-CNR and A-SCUC-CNR are the same.

The results presented in Table II and Table III prove that the proposed typical-decomposition and accelerated-decomposition methods are significantly faster for the same solution than extensive formulations of SCUC and SCUC-CNR respectively. It is intuitive that incorporating CNR will lead to additional computational complexity, which is demonstrated by the observation that the computing time of SCUC-CNR is longer than SCUC. However, it is the other way for the proposed approaches: the computational time for solving T/A-SCUC-CNR is much less than that for T/A-SCUC. The reason is that the addition of NR-PCFC sub-problem in addition to PCFC in T/A-SCUC-CNR leads to increased feasibility region of the sub-problems and reduced number of cuts and iterations. This is demonstrated and discussed in detail in sub-section VI.D.

TABLE II. SCUC ACCURACY ON THE IEEE 24-BUS SYSTEM

| MIPGAP=0.00 | SCUC | T-SCUC | A-SCUC |
|---|---|---|---|
| Total cost ($) | 963,893 | 963,893 | 963,893 |
| Solve time (s) | 6,013 | 2,440 | 1,351 |

TABLE III. SCUC-CNR ACCURACY ON THE IEEE 24-BUS SYSTEM

| MIPGAP=0.00 | SCUC-CNR | T-SCUC-CNR | A-SCUC-CNR |
|---|---|---|---|
| Total cost ($) | 928,794 | 928,794 | 928,794 |
| Solve time (s) | 9,625 | 47 | 9 |

### *B. AC Feasibility for CNR Solution*

Since the proposed methodologies are for a DC solution which is utilized in the industry for unit commitment, an AC feasibility check was performed to validate the CNR benefits. An example from the post-contingent sub-problems is presented in this subsection. Fig. 4 and Fig. 5 shows the MW flow as a percentage of line loading of lines close to the contingent line 27 in both DC solution and AC solution for IEEE 24-bus system, respectively. It was also noted that line 23, connecting the high-voltage and low-voltage regions, is the main bottle-neck line in the IEEE 24-bus system in majority of the critical contingencies. The CNR solution provided is line 37 as a congestion relief option in post-contingency phase. It can be observed from the line loading level for DC and AC solution that the benefits offered by CNR solution is not lost in an AC setting.

As stated earlier, the CNR action is only performed when there is a line congestion or overload after a contingency. From Fig. 4, we notice that the DC solution resulted in line overload in a post-contingency case for line 23 which is at 109%. The CNR action of switching line 37 resulted in flow redistribution in the network which relieves the line violation and brings the flow on line 23 to 100%.

An AC feasibility check for the above DC solution was verified and Fig. 5 represents the line loading in an AC setting. Here, line 23 is violated with a post-contingent flow of 108% after the contingency of line 27. Similarly, the CNR action of line 37 alleviates the violation on line 23 and results in a line-loading of 99.4%. This implies a similar decrease of ~9% is achieved with CNR action in both DC and AC setting for the violated line 23. The benefit of CNR is not lost in an AC solution as shown in Table IV. Moreover, the results showed that the MW flow on the lines were similar for most lines in both the AC and DC solutions. Therefore, the DC approximation holds well even when CNR is implemented.

TABLE IV. LINE LOADING OF LINE 23 IN IEEE 24-BUS SYSTEM

| Solution | Post-contingency | Post-reconfiguration |
|---|---|---|
| DC | 109 % | 100 % |
| AC | 108 % | 99.4 % |

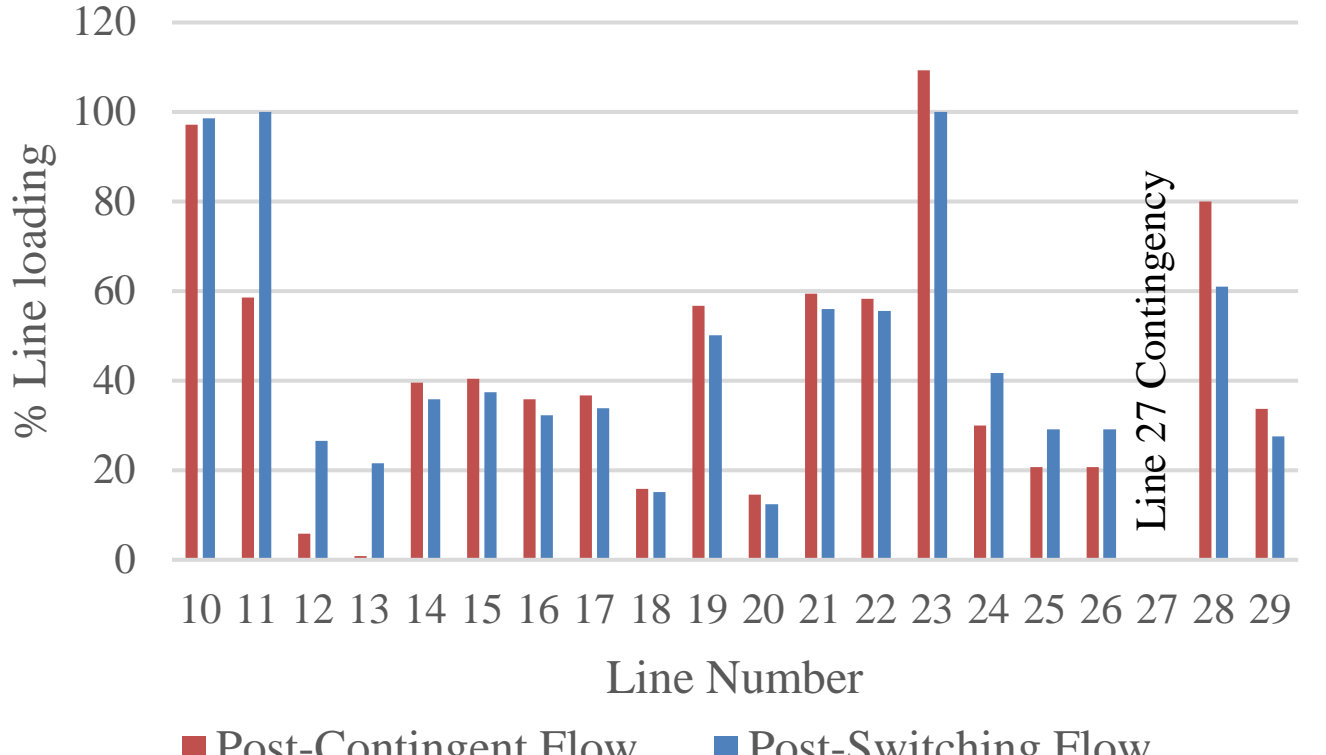

Fig. 4. Line loading for DC solution in IEEE 24-bus system.

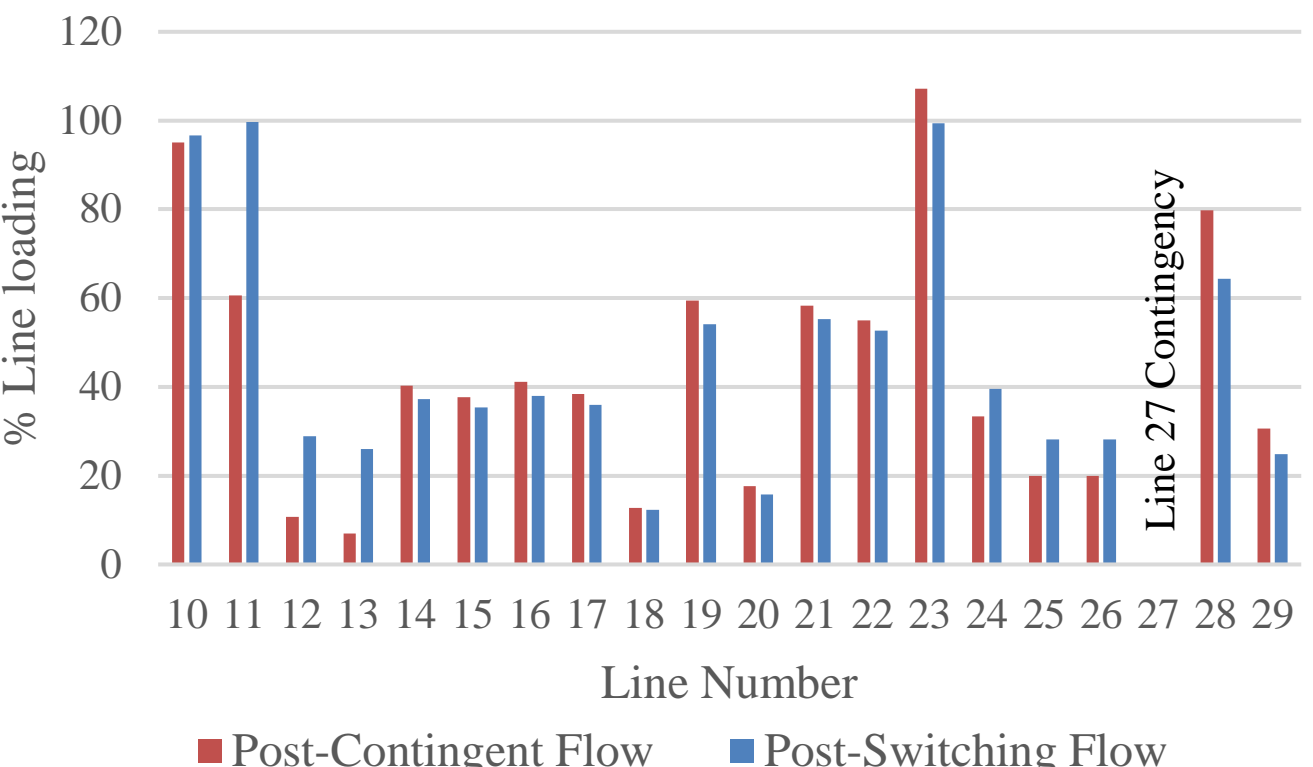

Fig. 5. Line loading for AC solution (MW only) in IEEE 24-bus system.

### *C. MIPGAP Sensitivity Analysis*

The MIPGAP, $\mu$, which is utilized for both the MUC and sub-problems affects the performance of all the methods. Typically, increasing the $\mu$ may lead to a less good solution, which increases the total cost. The impact of different selected $\mu$ values on total cost can be measured by the change in total cost in percentage ($\Delta Cost_{\mu}$), which is defined in (50).

$$\Delta Cost_{\mu} = \left( \frac{Cost_{\mu} - Cost_{\mu=0}}{Cost_{\mu=0}} \right) * 100\% \tag{50}$$

The MIPGAP sensitivity analysis is conducted on the IEEE 24-bus system and the results are shown in Fig. 6 and Fig. 7. Fig. 6 illustrates how $\Delta Cost_{\mu}$ varies with $\mu$ for different SCUC methods and SCUC-CNR methods. Fig. 7 shows that overall, the solve-time decreases significantly as $\mu$ increases. In general, extensive methods are more computationally intensive than the respective decomposition approaches except at very high relative MIPGAP such as $\mu$=1.0; higher optimality gap may result in a feasible solution for MUC faster but that solution can result in more violations in post-contingency sub-problem check which may result in additional cuts and iterations to converge. The extensive SCUC-CNR requires more computational time to execute than extensive SCUC. However, the opposite is true for the typical-decomposition and the accelerated-decomposition approach. This is because T/A-SCUC-CNR, (i) requires fewer iterations and (ii) the MUC size is marginally smaller compared to T/A-SCUC, as explained in sub-section IV.B. Another observation is that the solve time for T-SCUC-CNR and A-SCUC-CNR only reduces marginally as $\mu$ increases. One major reason is that sub-problems just verify post-contingency constraints for MUC commitment and dispatch schedule, and the solution always pointed to sub-problems being solved to $\mu$=0. Therefore, sub-problems are not affected by $\mu$. In addition, the MUC may only change marginally for the same test case.

Based on the above sensitivity analysis, $\mu$=0.01 (1%) provides a reasonable maximum cost gap of about 0.4% in a short time. For the rest of the paper, $\mu$=0.01 is used. However, it is worth noting that the performance of the proposed decomposition approaches implementing CNR fares well under tighter tolerances if higher accuracy is required.

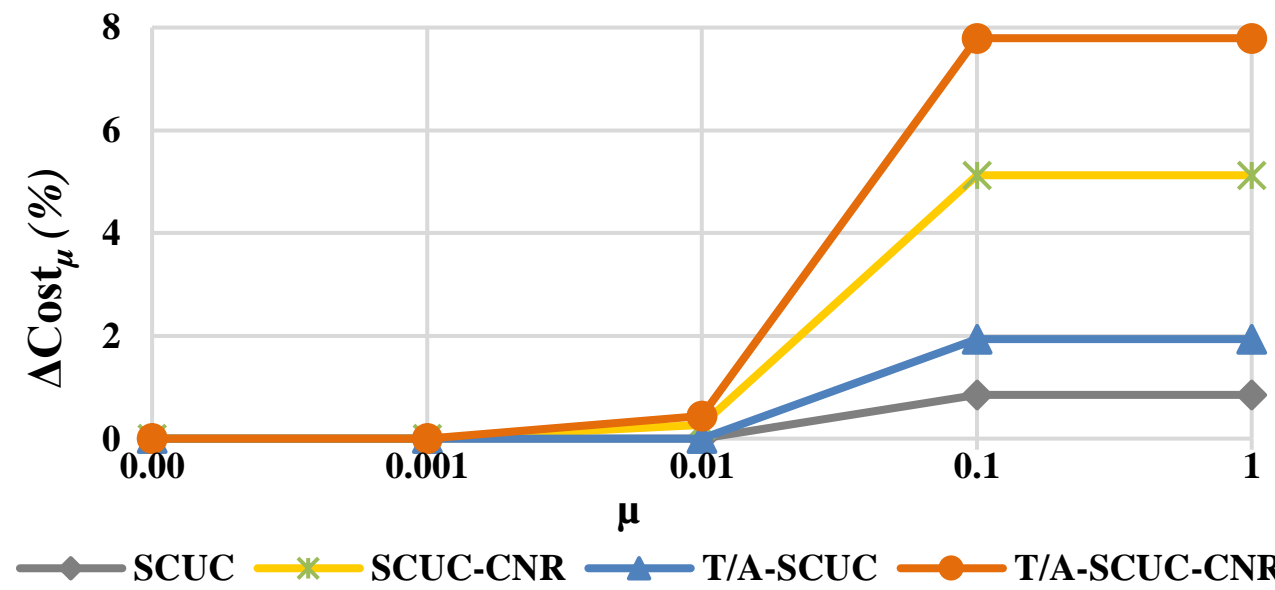

Fig. 6. $\Delta Cost_{\mu}$ versus relative MIPGAP $\mu$ on the IEEE 24-bus system.

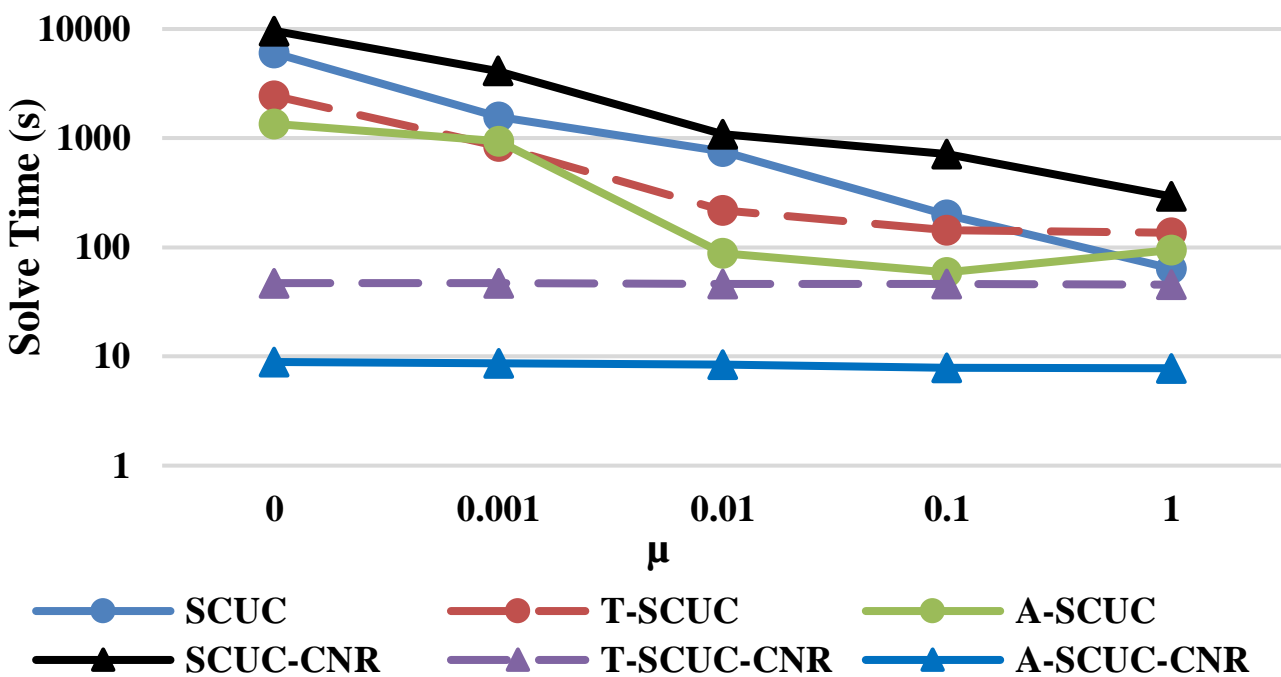

Fig. 7. Solve time versus relative MIPGAP $\mu$ on the IEEE 24-bus system.

### *D. Load Sensitivity Analysis*

Four scenarios were considered: two low-load/uncongested scenarios (80%, 90%), a base-load scenario (100%) and a high-load scenario (110%). The load profile was varied using a percentage multiplied to the nodal load. Table V shows the total cost for various methods under different load profiles for a MIPGAP=0.01

TABLE V. LOAD SENSITIVITY ANALYSIS ON IEEE 24-BUS SYSTEM

| Load Profile (%) | Total operational cost ($) | | | |
|---|---|---|---|---|
| | SCUC | T/A-SCUC | SCUC-CNR | T/A-SCUC-CNR |
| 80 | 467,883 | 467,883 | 467,883 | 467,883 |
| 90 | 624,398 | 624,398 | 623,459 | 623,459 |
| 100 | 963,893 | 963,893 | 932,919 | 932,919 |
| 110 | Infeasible | Infeasible | 1,424,140 | 1,424,140 |

In the low-load scenarios (80%, 90%), it is evident that CNR is never implemented as base-case network loading level is low and post-contingency networks are not congested. This implies all approaches obtain the same total cost.

CNR actions are observed in base-load and high-load scenarios (100%, 110%) where the network reconfiguration is utilized to relieve system congestion. This allows cheaper generators to produce more power, resulting in a reduced total operational cost. Interestingly, without CNR, the demand cannot be met due to network congestion for the high-load scenario.

## *E. Scalability Studies*

One of the key research gaps is the lack of an effective algorithm for solving SCUC-CNR that is scalable for large-scale power systems and solvable in realistic time. Table VI and Table VII tabulate the performance of SCUC and SCUC-CNR on IEEE 73-bus system respectively. Table V points that the extensive formulation of SCUC, requires a good starting point to solve in 7,743 seconds. One approach to have good starting solution is to utilize the commitment and dispatch results obtained from the relaxed MUC problem. However, without a starting solution, even SCUC proves to be infeasible in 100,000 seconds. A default starting solution can also be utilized where all generators are committed, which results in feasibility within 1% optimality gap in about 30,000 seconds that is still impractical. In the execution of the proposed decomposition approaches, a starting point solution is not considered and yet a feasible solution can be achieved faster. Based on Table VI, the starting point has a significant influence on a large optimization problem and it can be considered for T/A-SCUC and T/A-SCUC-CNR. Since the proposed decomposition approaches are iterative in nature, the best starting point can be obtained from the MUC solution from previous iteration. This may lead to further reduction in computational time. Also, the sub-problems are sequentially solved, and a parallel solving can speed up the algorithm.

Table VII shows that SCUC-CNR lacks scalability as it times out without a feasible solution for the IEEE 73-bus system when solved for 100,000 seconds with a good starting solution. However, this was bettered by T/A-SCUC-CNR.

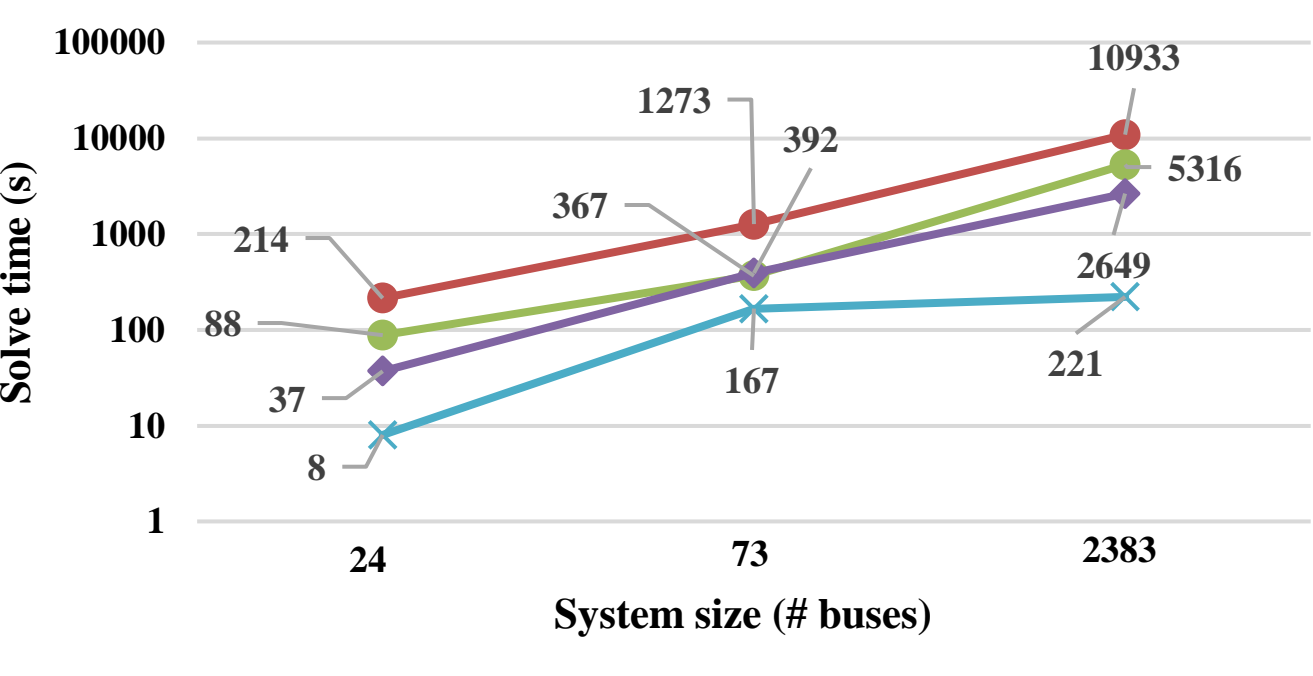

Fig. 8. Solving time versus system size.

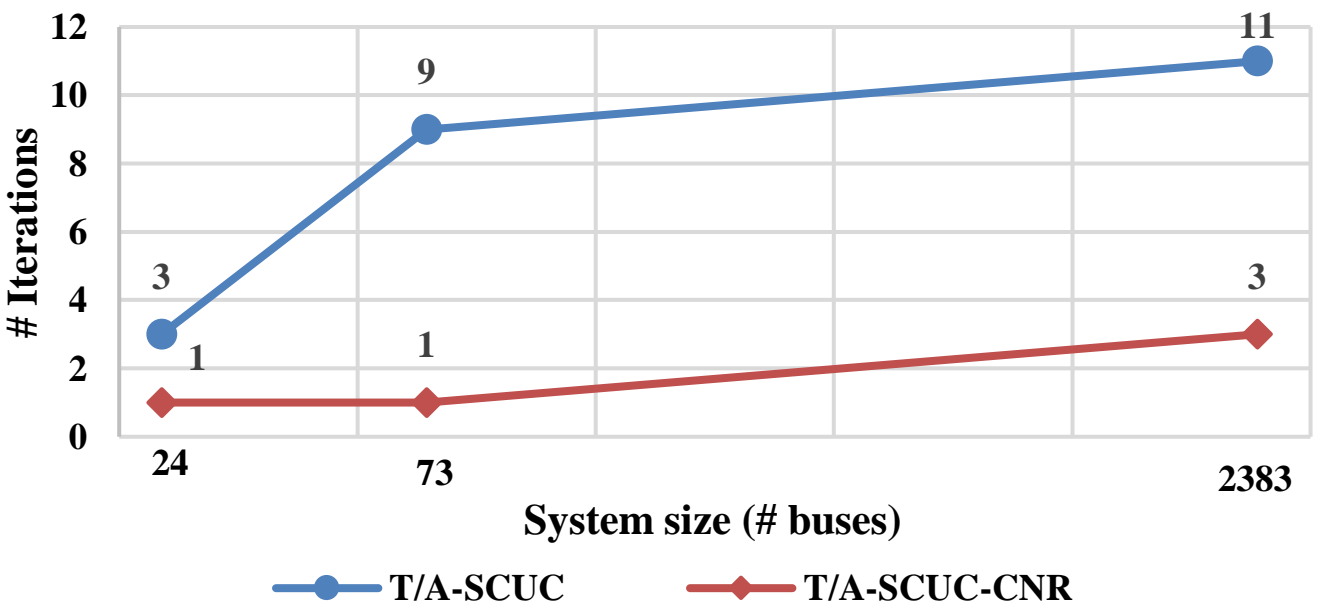

Fig. 9. Number of iterations versus size of the network.

TABLE VI. SCALABILITY OF SCUC TO IEEE 73-BUS SYSTEM

| MIPGAP=0.01 | SCUC | T-SCUC | A-SCUC |
|---|---|---|---|
| Total cost ($) | 3,224,980 | 3,223,760 | 3,223,760 |
| Solve time (s) | 7,743 | 1,273 | 367 |
| Feasibility | Feasible | Feasible | Feasible |
| Starting point | Yes | No | No |

TABLE VII. SCALABILITY SCUC-CNR TO IEEE 73-BUS SYSTEM

| MIPGAP=0.01 | SCUC-CNR | T-SCUC-CNR | A-SCUC-CNR |
|---|---|---|---|
| Total cost ($) | NA | 3,218,980 | 3,218,980 |
| Solve time (s) | 100,000 | 392 | 168 |
| Feasibility | Timeout | Feasible | Feasible |
| Starting point | Yes | No | No |

T/A-SCUC and T/A-SCUC-CNR are scalable to large networks such as the Polish system. Fig. 8 plots the solve time with respect to the size of the network. T/A-SCUC and T/A-SCUC-CNR are iterative in nature and Fig. 9 plots the number of iterations to solve the problem with respect to the size of the network. Due to the size of the Polish system, the 1-hour Polish case is utilized in Fig. 8 and Fig. 9 rather than the 24-hour Polish case to compare the performance with smaller systems. Here, it is noted that T/A-SCUC that do not perform CNR require more iterations to converge. The transmission flexibility obtained through implementing CNR in T/A-SCUC-CNR is evident from fewer iterations required to converge to a feasible solution with desired accuracy. This also means that the MUC problem that is more computationally intensive compared to sub-problems is solved fewer times, which saves a substantial amount of computational time. In addition, the number of cuts generated from infeasible post-contingency sub-problems for T/A-SCUC-CNR are also less than T/A-SCUC. In other words, the number of constraints added to the MUC problem for each iteration for T/A-SCUC-CNR is less than T/A-SCUC, which may lead to a less complex MUC problem and require less time to solve the MUC for each iteration for T/A-SCUC-CNR. The total number of cuts added to MUC for those decomposition methods is presented in Table VIII.

TABLE VIII. SUB-PROBLEM AND CUT DETAILS

| | IEEE 24-Bus | | IEEE 73-Bus | | Polish 1-Hr | |
|---|---|---|---|---|---|---|
| | T/A-SCUC | T/A-SCUC-CNR | T/A-SCUC | T/A-SCUC-CNR | T/A-SCUC | T/A-SCUC-CNR |
| # cuts | 198 | 42 | 65 | 17 | 76 | 14 |
| $\alpha$ | NA | 16 | NA | 20 | NA | 57 |

$\alpha$ in this table denotes number of sub-problems that were infeasible for post-contingency constraints without CNR but were feasible with CNR.

Table IX details the results of the Polish system when it is scaled to solve for 24-hour period. A-SCUC-CNR utilizes accelerators such as the CSPS and closest branches to contingency element (CBCE), a list of 20 closest lines to the contingent line. However, it was noted that over 96% off the 637 CNR actions, a top 10 choice from the CBCE list yields the solution. Therefore, the choice of 20 closest elements is a conservative approach. The inclusion of accelerators in A-SCUC-CNR decreases the solve time by 90% as compared to T-SCUC-CNR while the solution quality is retained. It is also evident that due to fewer iterations, T-SCUC-CNR is over 40% faster than A-SCUC-CNR.

TABLE IX. SCALABILITY TO POLISH SYSTEM FOR 24-HOUR PERIOD

| Parameters | A-SCUC | T-SCUC-CNR | A-SCUC-CNR |
|---|---|---|---|
| Total Cost ($) | 5,350,220 | 5,335,330 | 5,335,330 |
| € ($) | NA | 14,890 (0.28%) | 14,890 (0.28%) |
| Time (s) | 15,133.9 | 59,473.1 | 6,257.3 |
| $\delta$ | 0.04% | 0.12% | 0.12% |
| Iterations | 14 | 2 | 2 |
| # CNR | NA | 637 | 637 |
| # Cuts | 1,499 | 192 | 192 |

€ denotes the cost saving for T/A-SCUC-CNR as compared to A-SCUC. $\delta$ denotes the MIPGAP of the reported solution of MUC in the last iteration.

The comparison between A-SCUC and T/A-SCUC-CNR shows that there are 1,499 sub-problems resulting in cuts being added as constraints to the MUC problem for A-SCUC, as opposed to 192 cuts required for T/A-SCUC-CNR. Therefore, the MUC problem in A-SCUC is more constrained and takes longer to solve when compared to the MUC problem in T/A-SCUC-CNR. Not only that, the flexibility offered by CNR is evident by the following fact: out of 829 sub-problems that failed PCFC, 637 sub-problems are feasible with CNR through NR-PCFC, which implies about 77% of contingencies that failed feasibility in post-contingency check becomes feasible when CNR actions were implemented. Moreover, T/A-SCUC-CNR converge faster and require only 2 iterations against A-SCUC that requires 14 iterations. This implies the complex MUC problem is solved fewer times with T/A-SCUC-CNR leading to significant reduction in computational time.

The consideration of network reconfiguration for post-contingencies to alleviate network congestion in the large-scale Polish system for 24-hour period leads to a cost saving of $14,890. The discussion regarding the impact of CNR on congestion cost and markets is presented in detail in sub-section VI.F. It is to be noted that the results present an exhaustive monitoring of all non-radial transmission elements: 2,250 non-radial lines for the Polish system. This leads to 54,000 sub-problems for a 24-hour period per iteration whereas only 1,761 sub-problems were swiftly deemed as critical by CSPS. Subsequently, PCFC checked those 1,761 sub-problems and identified 829 sub-problems that failed feasibility check. The NR-PCFC that verifies feasibility of these contingencies with network reconfiguration further reduced the number of cuts required to be added to 192. Therefore, 637 sub-problems satisfied feasibility of post-contingency constraints by modifying the network topology. Though those 637 sub-problems that implemented CNR actions amounts to only 1.18% of all the sub-problems considered in the first iteration, considerable economic benefits are achieved with T/A-SCUC-CNR over A-SCUC. The solve time can be further significantly reduced if only a watch-list of key contingent lines are monitored as this will reduce the number of sub-problems drastically.

### *F. Comparison with State-of-the-art Method*

The state-of-the-art from [34], is based on nested-CCGA (NCCGA) and enhanced-nested-CCGA (E-NCCGA). Test results on the IEEE 118-bus system are utilized to compare the performance of the proposed algorithms against NCCGA and E-NCCGA. Here, it can be noted that the IEEE 118-bus system is easier to solve compared to the IEEE 73-bus system due to the smaller number of generators resulting in less binary variables and constraints for SCUC as seen in Table X. In addition, the IEEE 118-bus system has minimal generator data and a load profile leading to minimal system congestion.

TABLE X. EXTENSIVE METHOD PROBLEM SIZE

| Test System | IEEE 73-Bus | | IEEE 118-Bus | |
|---|---|---|---|---|
| Method | SCUC | SCUC-CNR | SCUC | SCUC-CNR |
| Constraints (C) | 1,555,824 | 2,517,192 | 1,477,176 | 3,716,212 |
| Equality C | 525,432 | 205,272 | 1,238,280 | 474,600 |
| Inequality C | 1,030,392 | 2,311,920 | 238,896 | 3,241,612 |
| Variables (V) | 799,515 | 1,121,763 | 1,318,512 | 2,066,630 |
| Binary V | 4,539 | 326,787 | 2,594 | 749,030 |
| Continuous V | 794,976 | 794,976 | 1,317,600 | 1,317,600 |
| Non-zeros | 4,716,249 | 7,914,081 | 4,505,304 | 11,955,218 |

Due to the simplistic nature of the above test system, an exact solution for both SCUC and SCUC-CNR was obtained easily. Based on the results shown in Table XI and Table XII, it is observed that T/A-SCUC and T/A-SCUC-CNR obtain the same exact solution, respectively. The results show that the IEEE 118-bus system only requires 165s (2.75 mins) and 59s (~1 min) for T-SCUC-CNR and A-SCUC-CNR, respectively. The T/A-SCUC-CNR is faster than the E-NCCGA and NCCGA which requires 2586s (43.1 mins) and 3504s (58.4 mins), respectively as reported in [34]. This points that the proposed method, A-SCUC-CNR is approximately 44 times faster than E-NCCGA and 59 times faster than NCCGA without compromising the optimal solution. However, it can be understood that T/A-SCUC-CNR only performs CNR and it might not be reasonable to compare against the E-NCCGA and NCCGA which performs both PNR and CNR for solve time. The modified IEEE 118-bus system and different hardware/software resources utilized in [34] can also affect the solution time.

TABLE XI. SCUC SOLUTIONS TO IEEE 118-BUS SYSTEM

| MIPGAP=0.01 | SCUC | T-SCUC | A-SCUC |
|---|---|---|---|
| Total cost ($) | 163,839 | 163,839 | 163,839 |
| Solve time (s) | 1,064 | 661 | 46 |
| Iterations | NA | 4 | 4 |

TABLE XII. SCUC-CNR SOLUTIONS TO IEEE 118-BUS SYSTEM

| MIPGAP=0.01 | SCUC-CNR | T-SCUC-CNR | A-SCUC-CNR |
|---|---|---|---|
| Total cost ($) | 163,609 | 163,609 | 163,609 |
| Solve time (s) | 3,045 | 165 | 59 |
| Iteration | NA | 1 | 1 |

### *G. Congestion Cost and Market Analysis*

The contingency-induced congestion cost, $CICC$, is calculated as the difference in total operation cost when emergency post-contingency line limits are imposed ($TC$) and not imposed ($TC_{NoEL}$) as represented in (51). The scenario when post-contingency emergency limits are not imposed is used as a benchmark since it is equivalent to implying that the system is not congested in the post-contingency situations. A-SCUC and A-SCUC-CNR are considered since we are interested in calculating the amount of $CICC$ reduced when CNR is implemented.

$$CICC = TC - TC_{NoEL} \tag{51}$$

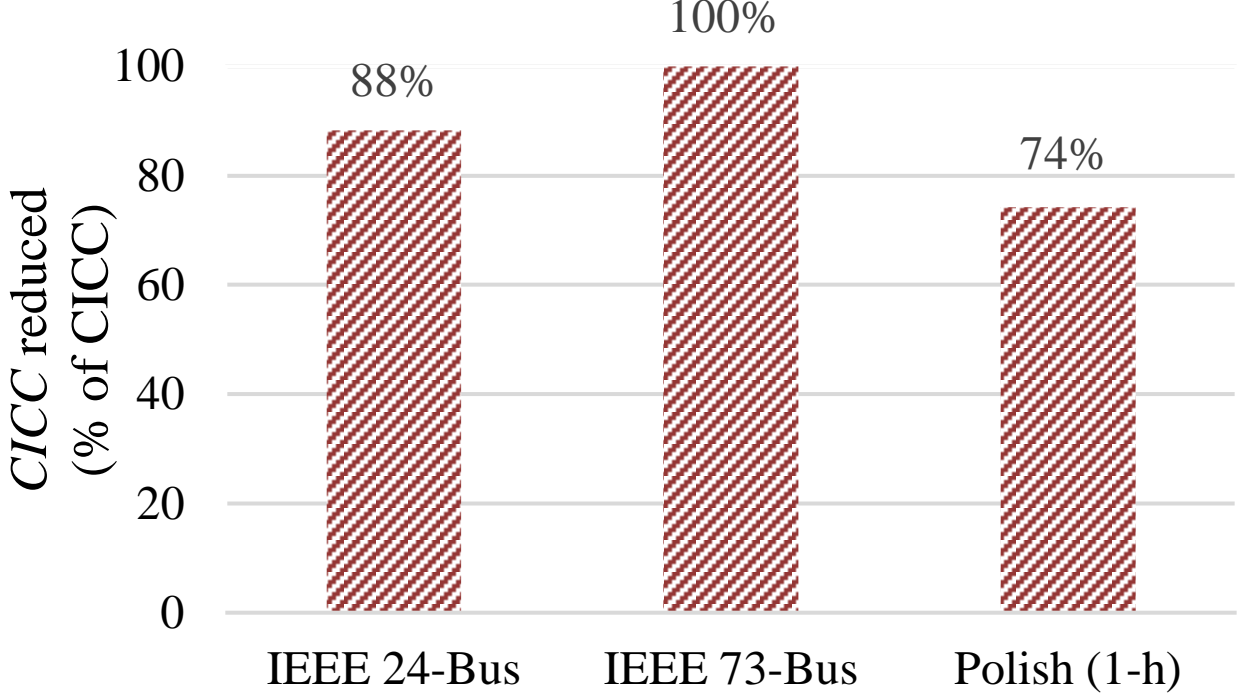

Fig. 10. *CICC* reduction with CNR for IEEE 24-bus, 73-Bus and Polish systems.

As observed from Fig. 10, the IEEE 24-bus system was the most congested system with a contingency-induced congestion cost of $35,099 due to the considered load profile along with lower transmission capability. This was followed by the 73-bus system and 1-hour polish system with $4,550 and $ 4,150 respectively. The *CICC* is considerably reduced with CNR in all the cases by 88%, 100% and 74% respectively. This is significant in heavily congested systems as seen in the case of IEEE 24-bus system where $30,974 is saved.

The market implication of reduction in *CICC* can be seen through the impact of CNR on nodal locational marginal prices (LMP). Table XIII shows the average nodal LMP calculated in various systems when CNR is not used and when CNR is implemented. Overall, it is observed that with CNR, (i) the average nodal LMP is reduced and (ii) the nodal LMP curve is flattened. It can be noted that congestion relief has a direct impact on the reduction in average nodal LMP. Similarly, it is also noted that the load payment is significantly reduced with CNR. Table XIV shows the total load payment for each test system with and without CNR. CNR resulted in a load payment reduction of $58,840 in the IEEE 24-bus system, $1,576,880 in the IEEE 73-bus system and $3,977 in the 1-hour Polish system, which correspond to percentage reductions of around 5.0%, 20.1% and 1.1% respectively. This makes sense since compared to the IEEE test systems, (i) the production cost of generators in the Polish system is low, (ii) the variation of system-wide generation cost in the Polish system is small, and (iii) the Polish system is loaded.

TABLE XIII. AVERAGE NODAL LMP ($/MWH)

| Test System | A-SCUC | | | | A-SCUC-CNR | | | |
|---|---|---|---|---|---|---|---|---|
| | Mean | Min | Max | StdD | Mean | Min | Max | StdD |
| IEEE 24-Bus | 23.39 | 5.46 | 150.6 | 0.86 | 23.23 | 5.46 | 150.6 | 0.84 |
| IEEE 73-Bus | 42.75 | 9.5 | 648.4 | 1.36 | 42.19 | 4.9 | 582.4 | 1.34 |
| Polish (1-hour) | 17.72 | 15.7 | 20.8 | 0.24 | 17.56 | 17.2 | 17.8 | 0.19 |

TABLE XIV. LOAD PAYMENT ($)

| Test System | A-SCUC | A-SCUC-CNR |
|---|---|---|
| IEEE 24 Bus | 1,171,220 | 1,112,380 |
| IEEE 73 Bus | 7,840,770 | 6,263,970 |
| Polish (1-hour) | 372,740 | 368,763 |

## VII. CONCLUSIONS

This paper proposes typical-decomposition and accelerated-decomposition approaches of SCUC and SCUC-CNR. The proposed decomposition approaches are generic and can be implemented to both SCUC and SCUC-CNR while outperforming the extensive formulations of SCUC and SCUC-CNR, respectively, in terms of (i) computational speed, (ii) algorithm scalability, and (iii) solution quality.

The accelerated-decomposition approach can easily link multiple accelerators to substantially reduce solution time. Specifically, CSPS, an exhaustive fast screening of sub-problems accurately identifies the critical contingent sub-problems which can lead to system overload or congestion. In addition, the A-SCUC-CNR benefits in computational speed achieved from ordered list, CBCE, for corrective actions.

The proposed A-SCUC-CNR utilizing the proposed accelerators, CSPS and CBCE, can solve a large-scale power system for 24-hour period in a reasonable time. As compared to T-SCUC-CNR, the proposed A-SCUC-CNR achieves a reduction of about 90% in the computational time without compromising solution accuracy. It can be noted that parallel solving of sub-problems was not considered in this paper and can provide additional computational time savings in future.

It was noted that implementation of CNR can achieve significant cost saving and provide feasible solutions for high critical demands where there are no feasible solutions without CNR. In addition, A-SCUC-CNR can provide quality solution much faster than the A-SCUC since fewer iterations are required. The load payment is dramatically reduced with CNR. Load payment reduction of 1%-20% can be realized for various networks. Mainly, the advantage of the proposed A-SCUC-CNR is that it provides quality solutions in a reasonable short time while dramatically reducing post-contingency network constraints induced congestion cost by 75%-100% in various scenarios. As a result, the total operation cost is reduced with CNR for congested networks.

## VIII. EXTENSIONS

### A. Forbidden Zones of Generators

Though the forbidden zones of generators are not considered in this work, they can be integrated easily in the iterative process. Some generators may consist of a few sub-regions between the minimum and maximum outputs where the generators are unstable; those unstable sub-regions are known as forbidden zones. For instance, as shown in Fig. 11, the sub-region between $P_g^{fmin}$ and $P_g^{fmax}$ is a forbidden zone. Typically, generators are allowed to operate in these zones provided the dispatch quickly moves away from this zone. If the generators cannot operate without crossing the region, then it is operated with full ramping rates inside the zone [37].

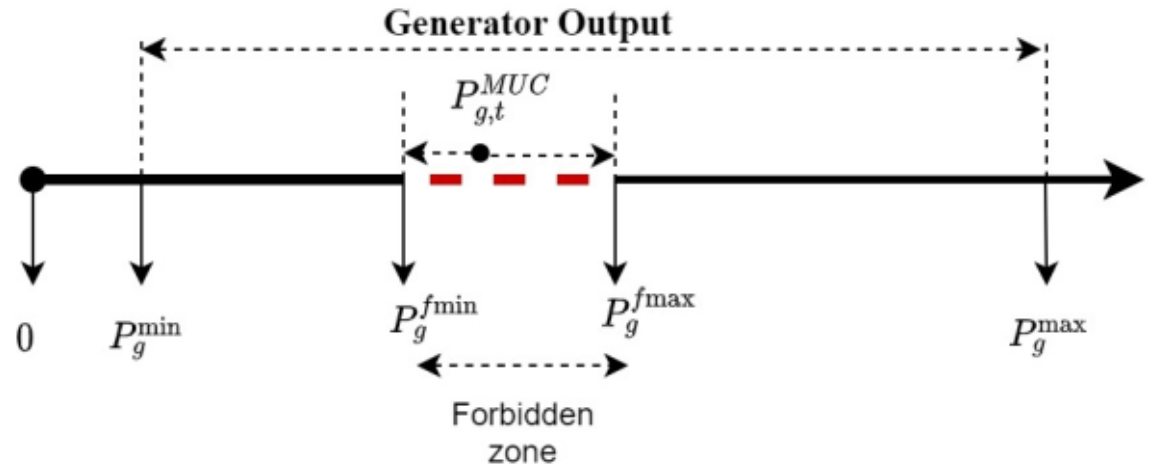

Fig. 11. Forbidden zones of operation in generators.

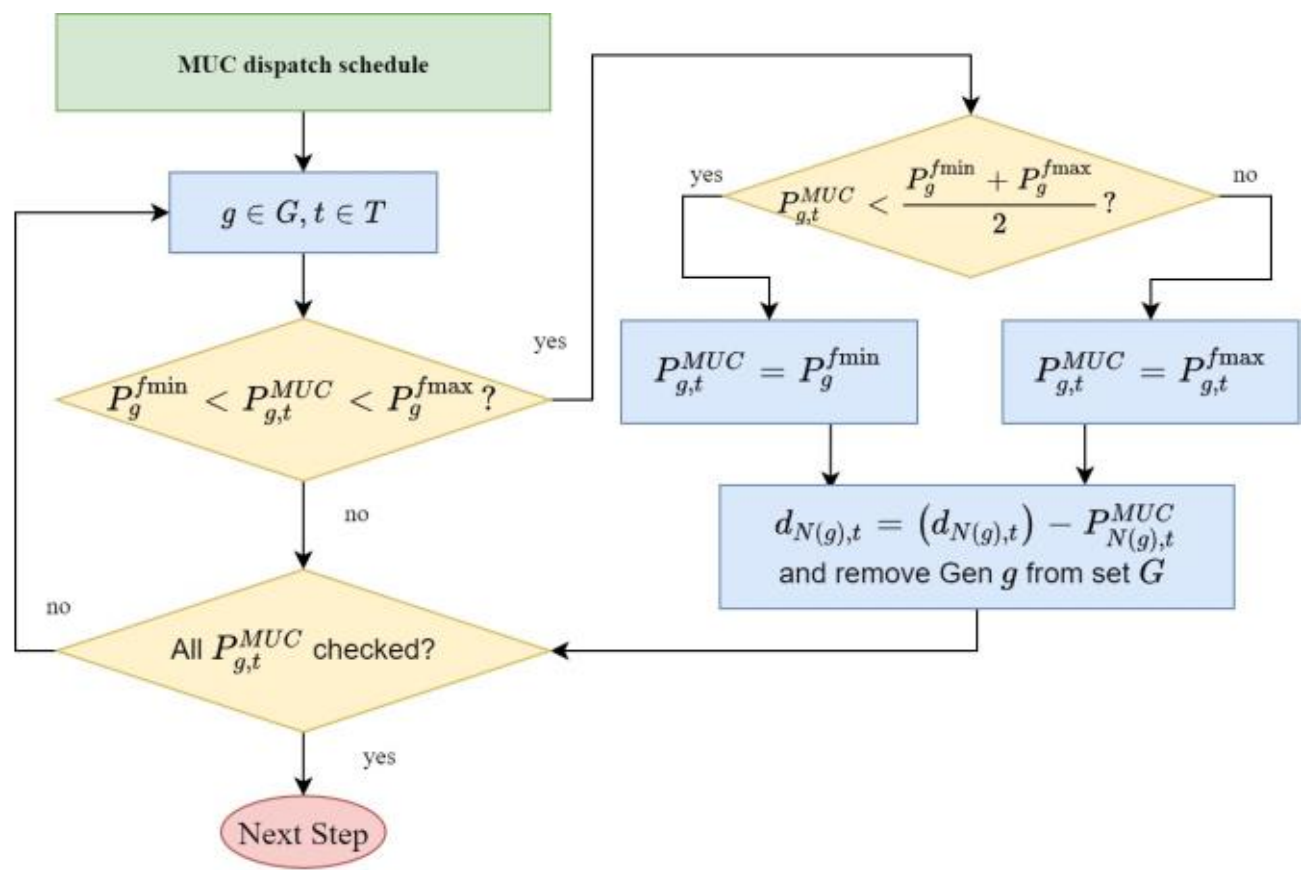

Fig. 12. Flowchart for considering forbidden zones in MUC.

Since the proposed method is an iterative process, the MUC solution after the final iteration can be verified against the forbidden zones for violations. The generator outputs for violated units can be modified as shown in Fig. 12. Once the violated generator outputs are fixed, the MUC is re-solved again for the rest of the units to satisfy the demand. Since the focus of this work is to investigate the application of CNR in SCUC, addressing generator forbidden zone is out of the scope of this paper.

### *B. Renewable Energy Source Integration*

The SCUC and SCUC-CNR models consider only hydro plants for this work. Other sources of variable renewable energy can be treated as traditional plants by considering a stochastic implementation which models multiple scenarios of weather as shown in [12]. This will convert deterministic SCUC/SCUC-CNR into stochastic SCUC/SCUC-CNR respectively. The commitment schedule obtained from stochastic SCUC/SCUC-CNR would be valid for all considered scenarios of weather. Since stochastic problems involve multiple scenarios, the resulting optimization problems consist of much more constraints than the respective deterministic problems. The complexity of such an extensive stochastic model can also be handled by decomposing it to a MUC model that provides the base-case commitment and dispatch solution and scenario-based post-contingency constraints that can be verified as sub-problems. Though it is evident that the consideration of scenarios will lead to additional sub-problems, it was showed that the proposed typical-decomposition and accelerated-decomposition methods can handle a large number of sub-problems in an effective manner. Note that the proposed decomposed approaches can be used to solve stochastic SCUC/SCUC-CNR with minor modifications, which will be our future work.

### *C. Other Future Work*

Similar to network reconfiguration, both demand response and energy storage can increase power system flexibility and can also be used as preventive or corrective actions [38]-[42]. Consideration of PNR and other technologies are considered as a future work. Though energy storage and demand response can be included explicitly by modelling constraints associated with their operational characteristics, they are not included in this paper since this paper focuses only on network reconfiguration as a corrective action. Prior research, [42]-[45], states that economic benefits can be obtained as a result of using energy storage and demand response in base-case. However, the effect of such technologies on post-contingent scenarios are not completely studied. Especially, the effectiveness of CNR versus corrective demand response and corrective energy storage actions can be studied and their interdependence can lead to more benefits, which will be a future scope considered as an extension of this paper.

CNR can also be implemented using bus-splitting with the formulations presented in [46]-[47]. This can map sub-station configuration using a node-breaker model and also consider bus-splitting contingencies which is another future extension.

**Arun Venkatesh Ramesh** (S’19) received the B.E. degree in electrical and electronic engineering from Anna University, Chennai, TN, India, and the M.Sc. degree in electrical engineering from Arizona State University, Tempe, AZ, USA, in 2013 and 2015 respectively.

Currently, he is pursuing his Ph.D. degree in electrical engineering at University of Houston, Houston, TX, USA. Prior to this, he worked with GE-Power, Hyderabad, India as a Software Engineer Specialist. His research interest includes power system operations and planning, numerical optimization, energy management system, energy markets, energy storage and renewable energy integration.

**Xingpeng Li** (S’13–M’18) received the B.S. degree in electrical engineering from Shandong University, Jinan, China, in 2010, and the M.S. degree in electrical engineering from Zhejiang University, Hangzhou, China, in 2013. He received the M.S. degree in industrial engineering and the Ph.D. degree in electrical engineering from Arizona State University, Tempe, AZ, USA, in 2016 and 2017 respectively.

Currently, he is an Assistant Professor in the Department of Electrical and Computer Engineering at the University of Houston. He previously worked for ISO New England, Holyoke, MA, USA, and PJM Interconnection, Audubon, PA, USA. Before joining the University of Houston, he was a Senior Application Engineer for ABB, San Jose, CA, USA. His research interests include power system operations, control and optimization, energy management system, energy markets, microgrids, and grid integration of renewable energy sources.

**Kory W. Hedman** (S’05–M’10–SM’16) received the B.S. degree in electrical engineering and the B.S. degree in economics from the University of Washington, the M.S. degree in economics and the M.S. degree in electrical engineering from Iowa State University, and the M.S. and Ph.D. degrees in operations research from the University of California at Berkeley.

He is currently the Director of the Power Systems Engineering Research Center, PSERC, and he previously served as a Program Director for the US Department of Energy (DOE) Advanced Research Projects Agency-Energy (ARPA-E). He was the recipient of the Presidential Early Career Award for Scientists and Engineers from U.S. President Barack H. Obama in 2017 and he is a prior recipient of the IEEE PES Outstanding Young Engineer Award. His research interests include power systems operations and planning, market design, power system economics, renewable energy, and operations research.